\documentclass[letterpaper, 10 pt, conference]{ieeeconf}  % Comment this line out if you need a4paper

\IEEEoverridecommandlockouts                              % This command is only needed if 
\usepackage{amsmath} % assumes amsmath package installed
\usepackage{amssymb}  % assumes amsmath package installed
\usepackage{mathtools}

\newcommand{\Span}[1]{\mathrm{span}\{#1\}}

\newcommand{\Dim}[1]{\mathrm{dim}(#1)}

\newcommand{\pfrac}[2]{\tfrac{\partial #1}{\partial #2}}
\newcommand{\ppfrac}[2]{\tfrac{\partial^2 #1}{\partial #2^2}}

\newcommand{\X}{\mathcal{X}}
\newcommand{\U}{\mathcal{U}}
\newcommand{\XU}{\mathcal{X}\times\mathcal{U}}

\newcounter{DefThmLemCorRem}
\newtheorem{definition}[DefThmLemCorRem]{Definition}
\newtheorem{theorem}[DefThmLemCorRem]{Theorem}

\newtheorem{proposition}[DefThmLemCorRem]{Proposition}
\newtheorem{remark}[DefThmLemCorRem]{Remark}

\title{\LARGE \bf
A Flat System Possessing no (x,u)-Flat Output
}

\author{Conrad Gst{\"o}ttner$^{1}$ Bernd Kolar$^{2}$ Markus Sch{\"o}berl$^{1}$% <-this % stops a space
\thanks{*The first author has been supported by the Austrian Science Fund (FWF) under grant number P 32151.}% <-this % stops a space
\thanks{$^{1}$Institute of Automatic Control and Control Systems Technology, Johannes Kepler University, Linz, Austria\newline {\tt\small \{conrad.gstoettner,markus.schoeberl\}@jku.at}}%
\thanks{$^{2}$Magna Powertrain Engineering Center Steyr GmbH \& Co KG, St. Valentin, Austria {\tt\small bernd\_kolar@ifac-mail.org}}%
}

\begin{document}

\maketitle
\thispagestyle{empty}
\pagestyle{empty}

%%%%%%%%%%%%%%%%%%%%%%%%%%%%%%%%%%%%%%%%%%%%%%%%%%%%%%%%%%%%%%%%%%%%%%%%%%%%%%%%
\begin{abstract}
	In general, flat outputs of a nonlinear system may depend on the system's state and input as well as on an arbitrary number of time derivatives of the latter. If a flat output which also depends on time derivatives of the input is known, one may pose the question whether there also exists a flat output which is independent of these time derivatives, i.e., an (x,u)-flat output. Until now, the question whether every flat system also possesses an (x,u)-flat output has been open. In this contribution, this conjecture is disproved by means of a counterexample. We present a two-input system which is differentially flat with a flat output depending on the state, the input and first-order time derivatives of the input, but which does not possess any (x,u)-flat output. The proof relies on the fact that every (x,u)-flat two-input system can be exactly linearized after an at most dim(x)-fold prolongation of one of its (new) inputs after a suitable input transformation has been applied. 
\end{abstract}

%%%%%%%%%%%%%%%%%%%%%%%%%%%%%%%%%%%%%%%%%%%%%%%%%%%%%%%%%%%%%%%%%%%%%%%%%%%%%%%%
\section{INTRODUCTION}

The concept of differential flatness has been introduced by Fliess, L\'evine, Martin and Rouchon in \cite{FliessLevineMartinRouchon:1992,FliessLevineMartinRouchon:1992-2}, and has attracted a lot of interest in the control systems community. The property of a system to be differentially flat (or just ``flat'' for short) allows for a systematic solution of feed-forward and feedback problems, see e.g. \cite{FliessLevineMartinRouchon:1995,DelaleauRudolph:1998,RudolphDelaleau:1998,FliessLevineMartinRouchon:1999}. Roughly speaking, a nonlinear control system of the form
	\begin{align}\label{eq:intro_nlsys}
		\begin{aligned}
			\dot{x}&=f(x,u)
		\end{aligned}
	\end{align}
	with $\Dim{x}=n$ states and $\Dim{u}=m$ inputs is flat if there exists an $m$-dimensional (fictitious) output 
	\begin{align}\label{eq:y}
		\begin{aligned}
			y&=\varphi(x,u,\dot{u},\ldots,u^{(q)})
		\end{aligned}
	\end{align}
	such that the state and the input of the system can locally be expressed as functions of this output and a finite number of its time derivatives, i.e.,
	\begin{align*}
		&\begin{aligned}
			x=F_x(y,\dot{y},\ldots,y^{(r-1)})
		\end{aligned}\\
		&\begin{aligned}
			u=F_u(y,\dot{y},\ldots,y^{(r)})\,.
		\end{aligned}
	\end{align*}
	Such a (fictitious) output \eqref{eq:y} is called a flat output of the system \eqref{eq:intro_nlsys}. The computation of flat outputs for systems of the general form \eqref{eq:intro_nlsys} is known to be a difficult problem and an active field of research since the introduction of flatness about 30 years ago. This problem, and flatness in general, have been studied within different mathematical frameworks. In \cite{FliessLevineMartinRouchon:1992}, \cite{FliessLevineMartinRouchon:1992-2} and \cite{FliessLevineMartinRouchon:1995}, a differential-algebraic setting is employed, whereas (infinite dimensional) differential-geometric frameworks including exterior differential systems are used e.g. in \cite{Pomet:1995,FliessLevineMartinRouchon:1999,MartinRouchon:1994,MartinRouchon:1995,MartinRouchon:1995-2,vanNieuwstadtRathinamMurray:1998,Kolar:2017}. Recent research in the field of flatness can be found e.g. in \cite{SchoberlSchlacher:2014,SilveiraPereiraRouchon:2015,LiNicolauRespondek:2016,GstottnerKolarSchoberl:2021-2,NicolauRespondek:2021}. Despite the use of sophisticated mathematical tools, up to now, there do not exist easily verifiable necessary and sufficient conditions for flatness, though, there are results available for several classes of systems, including systems which are linearizable by static feedback \cite{JakubczykRespondek:1980,HuntSu:1981}, two-input driftless systems \cite{MartinRouchon:1994}, systems linearizable by a one-fold prolongation of a suitably chosen input \cite{NicolauRespondek:2016,NicolauRespondek:2017}, two-input systems linearizable by a two-fold prolongation \cite{NicolauRespondek:2016-2,GstottnerKolarSchoberl:2021b}, and control-affine systems with four states and two inputs \cite{Pomet:1997}.
	
	One of the difficulties in checking flatness is that a flat output \eqref{eq:y} may a priori depend on time derivatives of $u$ up to an arbitrary order $q$. It is not known whether an upper bound on $q$ in terms of $n$ and/or $m$ exists.
	\begin{remark}
		Given a flat output $y=(\varphi_1,\ldots,\varphi_m)$, then for arbitrary $\beta>0$,
		\begin{align*}
			\begin{aligned}
				\bar{y}&=(\varphi_1,\varphi_2+\tfrac{\mathrm{d}^\beta}{\mathrm{d}t^\beta}\varphi_1,\varphi_3,\ldots,\varphi_m)
			\end{aligned}
		\end{align*}
		is also a flat output. By choosing $\beta$ large enough, we can thus easily construct a flat output which depends on time derivatives of $u$ up to an arbitrarily high order. Thus, a general bound on $q$ which covers all possible flat outputs of a system cannot exist. The question is whether every flat system admits a flat output $y=\varphi(x,u,\dot{u},\ldots,u^{(q)})$ with $q$ bounded in terms of $n$ and/or $m$.
	\end{remark}

	In \cite{MartinMurrayRouchon:2003}, it is conjectured that there might exist a bound on $q$ which is linear in the state dimension. This conjecture is motivated by the following example:
	\begin{align*}
		\begin{aligned}
			x_1^{(\alpha_1)}&=u_1
		\end{aligned}&
		\begin{aligned}
			x_2^{(\alpha_2)}&=u_2
		\end{aligned}
		\begin{aligned}
			\dot{x}_3&=u_1u_2
		\end{aligned}
	\end{align*}
	which admits the flat output
	\begin{align*}
		\begin{aligned}
			y_1&=x_3+\sum_{i=1}^{\alpha_1}(-1)^ix_1^{(\alpha_1-i)}u_2^{(i-1)}\\
			y_2&=x_2\,.
		\end{aligned}
	\end{align*}
	It is suspected that this system does not possess a flat output depending on derivatives of $u$ of order less than $\min(\alpha_1,\alpha_2)-1$, though, a proof of that has not been published. In fact, to the authors best knowledge, there has not even been published an example which is flat but evidentially does not admit any $(x,u)$-flat output (i.e., a flat output which may depend on the state $x$ and the input $u$, but not on time derivatives of $u$). In the following, we do exactly that. We prove that the system
	\begin{align}\label{eq:sys_xuup_intro}
		\begin{aligned}
			\dot{x}_1&=u_1\\
			\dot{x}_2&=u_2\\
			\dot{x}_3&=x_1+\tfrac{u_2^2}{2u_1}\,,
		\end{aligned}
	\end{align}
	which admits the $(x,u,\dot{u})$-flat output\footnote{The flat parameterization with respect to this flat output can be found in the appendix.}
	\vspace{-1ex}
	\begin{align}\label{eq:xuup_flat_output}
		\begin{aligned}
			y_1&=\tfrac{u_2}{u_1}\\
			y_2&=2x_2-\tfrac{2}{u_1}(x_1u_2-x_3\dot{u}_2)-2\tfrac{u_2}{u_1^2}(x_2\dot{u}_2+\\
			&\hspace{5em}x_3\dot{u}_1)+\tfrac{u_2^2}{u_1^3}(x_1\dot{u}_2+2x_2\dot{u}_1)-\tfrac{1}{u_1^4}x_1\dot{u}_1u_2^3
		\end{aligned}
	\end{align}
	cannot admit any $(x,u)$-flat output. Proofs in this paper are self contained, only some well known results about flatness are utilized without proving them.
\vspace{-0.5ex}
\section{PRELIMINARIES}
	In this section, we summarize some results regarding flatness of two-input systems. Throughout, all functions and vector fields are assumed to be smooth and all distributions are assumed to have locally constant dimension, we consider generic points only.
	
	Consider a nonlinear two-input system of the form
	\begin{align}\label{eq:sys}
		\begin{aligned}
			\dot{x}_i&=f_i(x,u)\,,&i=1,\ldots,n
		\end{aligned}
	\end{align}
	with the state $x$ being defined on an n-dimensional manifold $\X$ and with the input $u$ taking values in a two-dimensional manifold $\U$.
	\begin{definition}\label{def:flatness}
		The two-input system \eqref{eq:sys} is called flat if there exist two functions $y_j=\varphi_j(x,u,\dot{u},\ldots,u^{(q)})$, $j=1,2$ and smooth functions $F_{x_i}$ and $F_{u_j}$ such that locally
		\begin{align*}%\label{eq:param_r}
			\begin{aligned}
				x_i&=F_{x_i}(\varphi,\dot{\varphi},\ldots,\varphi^{(r-1)})\,,&i&=1,\ldots,n\\
				u_j&=F_{u_j}(\varphi,\dot{\varphi},\ldots,\varphi^{(r)})\,,&j&=1,2\,.
			\end{aligned}
		\end{align*}
		The functions $y_j=\varphi_j(x,u,\dot{u},\ldots,u^{(q)})$, $j=1,2$ are called the components of the $(x,u,\dot{u},\ldots,u^{(q)})$-flat output $y=\varphi(x,u,\dot{u},\ldots,u^{(q)})$.
	\end{definition}
	An important implication of Definition \ref{def:flatness} is that the time derivatives $\varphi,\dot{\varphi},\ldots,\varphi^{(\beta)}$ of a flat output are functionally independent for arbitrary $\beta$, which means that there does not exist any nontrivial function $\chi:\mathbb{R}^{2\beta+2}\mapsto\mathbb{R}$ such that $\chi(\varphi,\dot{\varphi},\ldots,\varphi^{(\beta)})=0$. (This property reflects the fact that there are no constraints on the time evolution of a flat output.) As a consequence of that, it can be shown that there exist unique minimal integers $r_1$,\,$r_2$ and unique maps $F_x$,\,$F_u$ such that
	\vspace{-1ex}
	\begin{align}\label{eq:param_R}
		\begin{aligned}
			x&=F_x(\varphi_1,\dot\varphi_1,\ldots,\varphi_1^{(r_1-1)},\varphi_2,\dot\varphi_2,\ldots,\varphi_2^{(r_2-1)})\\
			u&=F_u(\varphi_1,\dot\varphi_1,\ldots,\varphi_1^{(r_1)},\varphi_2,\dot\varphi_2,\ldots,\varphi_2^{(r_2)})\,.
		\end{aligned}
	\end{align}
	For $(x,u)$-flat outputs $y=\varphi(x,u)$, we define the integers $\rho_1$,\,$\rho_2$ where $\rho_j$ is the relative degree of the component $\varphi_j$ of the flat output, i.e., $\rho_j$ is the smallest integer such that $\varphi_j^{(\rho_j)}$ explicitly depends on $u$ (note that $\rho_j$ may be zero).
	
	Flat systems can be exactly linearized by means of an endogenous dynamic feedback, see e.g. \cite{FliessLevineMartinRouchon:1999} for the systematic construction of such a linearizing feedback. An important
	subclass of flat systems are those which are exactly linearizable by static feedback, i.e., which can be transformed into a linear controllable system, in particular the Brunovsk\'y normal form
	\vspace{-2ex}
	\begin{align*}
		\begin{aligned}
			\dot{\bar{x}}_{1,1}&=\bar{x}_{1,2}\\
			\dot{\bar{x}}_{1,2}&=\bar{x}_{1,3}\\
			&\vdotswithin{=}\\
			\dot{\bar{x}}_{1,\rho_1}&=\bar{u}_1
		\end{aligned}\qquad
		\begin{aligned}
			\dot{\bar{x}}_{2,1}&=\bar{x}_{2,2}\\
			\dot{\bar{x}}_{2,2}&=\bar{x}_{2,3}\\
			&\vdotswithin{=}\\
			\dot{\bar{x}}_{2,\rho_2}&=\bar{u}_2
		\end{aligned}
	\end{align*}
	where $\rho_1+\rho_2=n$, by means of an invertible state- and input transformation $\bar{x}=\Phi_x(x)$, $\bar{u}=\Phi_u(x,u)$. A system which is static feedback linearizable is obviously flat with $y=(\bar{x}_{1,1},\bar{x}_{2,1})$ as a flat output. The flat parameterization \eqref{eq:param_R} with respect to this flat output is a diffeomorphism and conversely, if a system possesses a flat output for which the corresponding flat parameterization is a diffeomorphism, then it is static feedback linearizable. We refer to these particular flat outputs as linearizing outputs. The relative degrees $\rho_1$,\,$\rho_2$ of the components of linearizing outputs sum to $n$ (a linearizing output is an output with a vector relative degree of $n$, see e.g. \cite{Isidori:1995}). The static feedback linearization problem has been solved completely, see \cite{JakubczykRespondek:1980}, \cite{HuntSu:1981} and \cite{NijmeijervanderSchaft:1990}, \cite{Isidori:1995}. In the following we recall necessary and sufficient conditions for static feedback linearizability of systems of the form \eqref{eq:sys}, for a proof we refer to \cite{NijmeijervanderSchaft:1990}.
	For \eqref{eq:sys}, define the distributions $\mathcal{D}^0=\Span{\partial_{u_1},\partial_{u_2}}$ and $\mathcal{D}^i=\mathcal{D}^{i-1}+[f,\mathcal{D}^{i-1}]$, $i\geq 1$ on the state and input manifold $\XU$, where $f=f_i(x,u)\partial_{x_i}$. 
	\begin{theorem}\label{thm:sfl}
		The two-input system \eqref{eq:sys} is linearizable by static feedback if and only if all the distributions $\mathcal{D}^i$ are involutive and $\Dim{\mathcal{D}^n}=n+2$.
	\end{theorem}
	The following proposition states that two-input $(x,u)$-flat systems can be exactly linearized by a special kind of endogenous dynamic feedback, namely prolongations of inputs. Based on this result, we will prove that \eqref{eq:sys_xuup_intro} cannot admit any $(x,u)$-flat output.
	\begin{proposition}\label{prop:linearization_xu}
		Every $(x,u)$-flat two-input system \eqref{eq:sys} can be rendered static feedback linearizable by $\Dim{x}=n$-fold prolonging a suitably chosen input after a suitable invertible input transformation has been applied.
	\end{proposition}
	\textit{Proof.} Let $y=\varphi(x,u)$ be an $(x,u)$-flat output of \eqref{eq:sys} and let $r_1$,\,$r_2$ be the unique smallest integers such that \eqref{eq:param_R} holds. The derivatives $\varphi_1,\dot\varphi_1,\ldots,\varphi_1^{(\rho_1-1)},\varphi_1,\dot\varphi_2,\ldots,\varphi_2^{(\rho_2-1)}$ yield in total $\rho_1+\rho_2$ independent functions of $x$ only. Apply the invertible input transformation $\bar u_1=\varphi_1^{(\rho_1)}(x,u)$, $\bar u_2=u_2$ (permute if necessary $u_1$ and $u_2$). %, and permute furthermore $\varphi_1$ and $\varphi_2$ if $\pfrac{\varphi_1}{u}(x_0,u_0)=0$).
	After applying this transformation, we clearly have $\varphi_2^{(\rho_2)}=\varphi_2^{(\rho_2)}(x,\bar{u}_1)$. Otherwise, the flat output and its derivatives up to arbitrary orders would only yield $\rho_1+\rho_2$ independent functions of $x$ and we could not express all states in terms of the flat output and its derivatives, which contradicts flatness. (An exception is of course $\varphi$ being a linearizing output, i.e., $\rho_1+\rho_2=n$, in which case the proposition clearly holds, since $n$-fold prolonging an input of a system in Brunovsk\'y normal form again yields a system in Brunovsk\'y normal form.) So we have $\varphi_2^{(\rho_2)}=\varphi_2^{(\rho_2)}(x,\bar{u}_1)$, and it follows that
	\begin{align}\label{eq:phi_R}
		\begin{aligned}
			\varphi_1&=\varphi_1(x)&&&\varphi_2&=\varphi_2(x)\\
			&\vdotswithin{=}&&&&\vdotswithin{=}\\
			\varphi_1^{(\rho_1-1)}&=\varphi_1^{(\rho_1-1)}(x)&&&\varphi_2^{(\rho_2-1)}&=\varphi_2^{(\rho_2-1)}(x)\\
			\varphi_1^{(\rho_1)}&=\bar u_1&&&\varphi_2^{(\rho_2)}&=\varphi_2^{(\rho_2)}(x,\bar u_1)\\
			\varphi_1^{(\rho_1+1)}&=\dot{\bar{u}}_1&&&\varphi_2^{(\rho_2+1)}&=\varphi_2^{(\rho_2+1)}(x,\bar u_1,\dot{\bar{u}}_1)\\
			&\vdotswithin{=}&&&&\vdotswithin{=}\\
			\varphi_1^{(r_1)}&=\bar u_1^{(r_1-\rho_1)}&&&\varphi_2^{(r_2)}&=\varphi_2^{(r_2)}(x,\bar u_1,\dot{\bar{u}}_1,\ldots,\\
			&&&&&\hspace{5em}\bar{u}_1^{(r_2-\rho_2)},\bar{u}_2)
		\end{aligned}
	\end{align}
	with $\varphi_2^{(\rho_2+k)}$ explicitly depending on $\bar u_1^{(k)}$. To show that $\varphi_2^{(r_2)}$ explicitly depends on $\bar{u}_2$ but no lower order time derivative of $\varphi_2$ can depend on $\bar{u}_2$, recall that $r_2$ is by assumption the minimal integer for which \eqref{eq:param_R} holds. If $\bar{u}_2$ would occur in $\varphi_2^{(s)}$ for some $s<r_2$, then $r_2$ could not be minimal since $\varphi_2^{(s+1)},\ldots,\varphi_2^{(r_2)}$ would be useless for constructing functions of $x$ and $\bar{u}$ only. On the other hand, $\bar{u}_2$ must explicitly occur in $\varphi_2^{(r_2)}$, since otherwise $\bar{u}_2$ could not be expressed in terms of $\varphi_1,\ldots,\varphi_1^{(r_1)},\varphi_2,\ldots,\varphi_2^{(r_2)}$. Furthermore, it follows that the explicit dependence of $\varphi_2^{(r_2)}$ on $\bar{u}_1^{(r_2-\rho_2)}$ and the minimality of $r_1$ imply that $r_1-\rho_1=r_2-\rho_2$. From the functions \eqref{eq:phi_R}, we obtain exactly $\rho_1+r_2$ independent functions of $x$ only, and by the flatness assumption this must be equal to $n$, i.e., $\rho_1+r_2=n$. Therefore, the $r_1+r_2+2$ functions \eqref{eq:phi_R} depend on the $n+2+(r_1-\rho_1)=r_1+r_2+2$ variables $x,\bar{u}_1,\dot{\bar{u}}_1,\ldots,\bar{u}_1^{(r_1-\rho_1)},\bar{u}_2$, and since the time derivatives of a flat output up to an arbitrary order are functionally independent, it follows that conversely the variables $x,\bar{u}_1,\dot{\bar{u}}_1,\ldots,\bar{u}_1^{(r_1-\rho_1)},\bar{u}_2$ can be expressed in terms of the functions \eqref{eq:phi_R}. In other words, \eqref{eq:phi_R} describes a diffeomorphism. Above, we noted that $r_1-\rho_1=r_2-\rho_2$ and $\rho_1+r_2=n$. We thus also have $\rho_2+r_1=n$, and since $\rho_1,\rho_2\geq 0$, we clearly have $r_1,r_2\leq n$ and in particular $r_1-\rho_1\leq n$. Therefore, we may extend \eqref{eq:phi_R} by adding
	\begin{align}\label{eq:extension}
		\begin{aligned}
			\varphi_1^{(r_1+1)}&=\bar{u}_1^{(r_1-\rho_1+1)}\\
			&\vdotswithin{=}\\
			\varphi_1^{(n+\rho_1)}&=\bar{u}_1^{(n)}
		\end{aligned}
	\end{align}
	and still have a diffeomorphism (note that in case of $r_1-\rho_1=n$ there are no equations added). Now consider the prolonged system
	\vspace{-1ex}
	\begin{align}\label{eq:sys_p}
		\begin{aligned}
			\dot{x}&=f(x,\bar{u}_1,\bar{u}_2)\\
			\dot{\bar{u}}_1&=\bar{u}_{1,1}\\
			\dot{\bar{u}}_{1,1}&=\bar{u}_{1,2}\\
			&\vdotswithin{=}\\
			\dot{\bar{u}}_{1,n-1}&=\bar{u}_{1,n}\\
		\end{aligned}
	\end{align}
	with the state $(x,\bar{u}_1,\bar{u}_{1,1},\ldots,\bar{u}_{1,n-1})$ and the input $(\bar{u}_{1,n},\bar{u}_2)$. The flat parameterization of \eqref{eq:sys_p} with respect to the flat output $\varphi$ is a diffeomorphism (the inverse of which is given by \eqref{eq:phi_R} extended by \eqref{eq:extension} with $\bar{u}_1^{(\alpha)}$ replaced by $\bar{u}_{1,\alpha}$), and thus \eqref{eq:sys_p} it is static feedback linearizable with $\varphi$ as a linearizing output.\hfill$\square$
\begin{remark}
	In \cite{GstottnerKolarSchoberl:2020} it has been shown that every two-input $(x,u)$-flat system \eqref{eq:sys} can be rendered static feedback linearizable by an $(r_1+r_2-n)$-fold prolongation of a suitably chosen input after a suitable input transformation has been applied. In \cite{KolarSchoberlSchlacher:2016} bounds on $r_j$ in terms of the integer $q$ in $y=\varphi(x,u,\dot{u},\ldots,u^{(q)})$ have been derived. For $(x,u)$-flat outputs, i.e., $q=0$, this bound is given by $r_j\leq n$ and hence $r_1+r_2-n\leq n$. These two results thus imply that at most an $n$-fold prolongation is required for rendering a two-input $(x,u)$-flat system static feedback linearizable. However, Proposition \ref{prop:linearization_xu} is not directly implied by these two results since Proposition \ref{prop:linearization_xu} states that a linearization is always possible with exactly $n$ prolongations (even though the minimal number of prolongations required for the considered system may be less than $n$).
\end{remark}

\section{EXAMPLE}
	Consider again the system \eqref{eq:sys_xuup_intro}
%	\begin{align}\label{eq:sys_xuup_intro}
%		\begin{aligned}
%			\dot{x}_1&=u_1\\
%			\dot{x}_2&=u_2\\
%			\dot{x}_3&=x_1+\tfrac{u_2^2}{2u_1}
%		\end{aligned}
%	\end{align}
	from the introduction. This system admits an $(x,u,\dot{u})$-flat output (given in the introduction), and it can be shown that the system is linearizable by a $4$-fold prolongation of the input $\bar{u}_1=\tfrac{u_2}{u_1}$ (see appendix). In the following, based on Proposition \ref{prop:linearization_xu}, we show that the system cannot be $(x,u)$-flat.
%	We claim that the system does not possess an $(x,u)$-flat output.
	\begin{proposition}
		The system \eqref{eq:sys_xuup_intro} does not possess any $(x,u)$-flat output.
	\end{proposition}
	\textit{Proof.} We prove this proposition by contradiction. Assume that \eqref{eq:sys_xuup_intro} is $(x,u)$-flat. Then, according to Proposition \ref{prop:linearization_xu}, the system can be rendered static feedback linearizable by $n=3$-fold prolonging a suitable input after a suitable invertible input transformation has been applied. The most general form of such an input transformation is 
	\begin{align}\label{eq:input_trasnformation}
		\begin{aligned}
			\bar{u}_1&=g_1(x,u_1,u_2)\\
			\bar{u}_2&=g_2(x,u_1,u_2)\,,
		\end{aligned}
	\end{align}
	with an inverse of the general form
	\begin{align}\label{eq:input_trasnformation_inv}
		\begin{aligned}
			u_1&=h_1(x,\bar{u}_1,\bar{u}_2)\\
			u_2&=h_2(x,\bar{u}_1,\bar{u}_2)\,.
		\end{aligned}
	\end{align}
	By assumption, there exists a transformation of the form \eqref{eq:input_trasnformation} with inverse \eqref{eq:input_trasnformation_inv} such that the prolonged system
	\begin{align*}
		\begin{aligned}
			\dot{x}_1&=h_1(x,\bar{u}_1,\bar{u}_2)\\
			\dot{x}_2&=h_2(x,\bar{u}_1,\bar{u}_2)\\
			\dot{x}_3&=x_1+\tfrac{(h_2(x,\bar{u}_1,\bar{u}_2))^2}{2h_1(x,\bar{u}_1,\bar{u}_2)}\\
			\dot{\bar{u}}_1&=\bar{u}_{1,1}\\
			\dot{\bar{u}}_{1,1}&=\bar{u}_{1,2}\\
			\dot{\bar{u}}_{1,2}&=\bar{u}_{1,3}
		\end{aligned}
	\end{align*}
	is static feedback linearizable. By the regularity of \eqref{eq:input_trasnformation_inv}, we always have at least $\pfrac{h_1}{\bar{u}_2}\neq 0$ or $\pfrac{h_2}{\bar{u}_2}\neq 0$. This allows for at least one of the normalizations
	\begin{align}\label{eq:prolonged_systems}
		\begin{aligned}
			\dot{x}_1&=\bar{u}_2\\
			\dot{x}_2&=h(x,\bar{u}_1,\bar{u}_2)\\
			\dot{x}_3&=x_1+\tfrac{(h(x,\bar{u}_1,\bar{u}_2))^2}{2\bar{u}_2}\\
			\dot{\bar{u}}_1&=\bar{u}_{1,1}\\
			\dot{\bar{u}}_{1,1}&=\bar{u}_{1,2}\\
			\dot{\bar{u}}_{1,2}&=\bar{u}_{1,3}
		\end{aligned}&&\text{or}&&\begin{aligned}
			\dot{x}_1&=h(x,\bar{u}_1,\bar{u}_2)\\
			\dot{x}_2&=\bar{u}_2\\
			\dot{x}_3&=x_1+\tfrac{\bar{u}_2^2}{2h(x,\bar{u}_1,\bar{u}_2)}\\
			\dot{\bar{u}}_1&=\bar{u}_{1,1}\\
			\dot{\bar{u}}_{1,1}&=\bar{u}_{1,2}\\
			\dot{\bar{u}}_{1,2}&=\bar{u}_{1,3}\,,
		\end{aligned}
	\end{align}
	where by abuse of notation we denote the newly introduced input still by $\bar{u}_2$ and renamed the composition of the remaining function ($h_2$ or $h_1$) with the ``inverse'' of this normalization by $h$ (without a subscript). This normalization can be done before prolonging the input $\bar{u}_1$ since it does not involve time derivatives of $\bar{u}_1$.
	
	In the following, we show that there cannot exist an invertible input transformation \eqref{eq:input_trasnformation_inv} such that at least one of the prolonged systems \eqref{eq:prolonged_systems} is static feedback linearizable. This will allow us to conclude that the system is not linearizable by a three-fold prolongation and consequently not $(x,u)$-flat.\\
	
	\textit{\textbf{Case 1.}} Let us first consider the prolonged system
	\begin{align}\label{eq:prolonged_1}
		\begin{aligned}
			\dot{x}_1&=\bar{u}_2\\
			\dot{x}_2&=h(x,\bar{u}_1,\bar{u}_2)\\
			\dot{x}_3&=x_1+\tfrac{(h(x,\bar{u}_1,\bar{u}_2))^2}{2\bar{u}_2}\\
			\dot{\bar{u}}_1&=\bar{u}_{1,1}\\
			\dot{\bar{u}}_{1,1}&=\bar{u}_{1,2}\\
			\dot{\bar{u}}_{1,2}&=\bar{u}_{1,3}\,,
		\end{aligned}
	\end{align}
	where the corresponding input transformation reads as	
	\begin{align}\label{eq:inv_trans_1}
		\begin{aligned}
			u_1&=\bar{u}_2\\
			u_2&=h(x,\bar{u}_1,\bar{u}_2)\,.
		\end{aligned}
	\end{align}
	By assumption, the prolonged system \eqref{eq:prolonged_1} is static feedback linearizable. Thus, according to Theorem \ref{thm:sfl}, the distributions
	\begin{align*}
		\begin{aligned}
			\mathcal D_p^0&=\Span{\partial_{\bar{u}_{1,3}},\partial_{\bar{u}_2}}\\
			\mathcal D_p^1&=\Span{\partial_{\bar{u}_{1,3}},\partial_{\bar{u}_{1,2}},\partial_{\bar{u}_2},\partial_{x_1}\!\!+\!\pfrac{h}{\bar{u}_2}\partial_{x_2}\!\!+\!\tfrac{2h\pfrac{h}{\bar{u}_2}\bar{u}_2-h^2}{2\bar{u}_2^2}\partial_{x_3}}\\
			&\vdotswithin{=}
		\end{aligned}
	\end{align*}
	are involutive. The latter implies that
	\begin{align*}
		[\partial_{\bar{u}_2},\partial_{x_1}+\pfrac{h}{\bar{u}_2}\partial_{x_2}+\tfrac{2h\pfrac{h}{\bar{u}_2}\bar{u}_2-h^2}{2\bar{u}_2^2}\partial_{x_3}]\in\mathcal D_p^1\,,
	\end{align*}
	which can only hold if 
	\begin{align*}
		\begin{aligned}
			\ppfrac{h}{\bar{u}_2}&=0&\text{and}&&\pfrac{}{\bar{u}_2}(\tfrac{2h\pfrac{h}{\bar{u}_2}\bar{u}_2-h^2}{2\bar{u}_2^2})&=0
		\end{aligned}
	\end{align*}
	The first condition implies that $h$ is affine with respect to $\bar{u}_2$, i.e., we actually have $h=a(x,\bar{u}_1)+b(x,\bar{u}_1)\bar{u}_2$. The second condition then yields
	\begin{align*}
		\begin{aligned}
			\pfrac{}{\bar{u}_2}\left(\frac{b^2\bar{u}_2^2-a^2}{2\bar{u}_2^2}\right)&=0
		\end{aligned}
	\end{align*}
	which simplifies to $\tfrac{a^2}{\bar{u}_2^3}=0$ and thus $a=0$. The distribution $\mathcal D_p^1$ thus simplifies to
	\begin{align*}
		\begin{aligned}
			\mathcal D_p^1&=\Span{\partial_{\bar{u}_{1,3}},\partial_{\bar{u}_{1,2}},\partial_{\bar{u}_2},\partial_{x_1}+b\partial_{x_2}+\tfrac{b^2}{2}\partial_{x_3}}\,,
		 \end{aligned}
 	\end{align*}
 	and for the next distribution we obtain
	\begin{align*}
	 	\begin{aligned}
		 	\mathcal D_p^2&=\Span{\partial_{\bar{u}_{1,3}},\partial_{\bar{u}_{1,2}},\partial_{\bar{u}_{1,1}},\partial_{\bar{u}_2},\partial_{x_1}+b\partial_{x_2}+\tfrac{b^2}{2}\partial_{x_3},\\
		 		&\hspace{0.4ex}(x_1\pfrac{b}{x_3}+\bar{u}_{1,1}\pfrac{b}{\bar{u}_1})\partial_{x_2}+(x_1b\pfrac{b}{x_3}+\bar{u}_{1,1}b\pfrac{b}{\bar{u}_1}-1)\partial_{x_3}}\,,
	 	\end{aligned}
 	\end{align*}
 	which of course must again be involutive. The latter in particular implies that
 	\begin{align*}
 		&[\partial_{\bar{u}_{1,1}},(x_1\pfrac{b}{x_3}+\bar{u}_{1,1}\pfrac{b}{\bar{u}_1})\partial_{x_2}+\\
 		&\hspace{5em}(x_1b\pfrac{b}{x_3}+\bar{u}_{1,1}b\pfrac{b}{\bar{u}_1}-1)\partial_{x_3}]\in \mathcal D_p^2\,,
 	\end{align*}
 	i.e., $\pfrac{b}{\bar{u}_1}\partial_{x_2}+b\pfrac{b}{\bar{u}_1}\partial_{x_3}\in \mathcal D_p^2$, which can only hold if the vector field
 	\begin{align*}
 		\begin{aligned}
 			\pfrac{b}{\bar{u}_1}\partial_{x_2}+b\pfrac{b}{\bar{u}_1}\partial_{x_3}
 		\end{aligned}
 	\end{align*}
 	is collinear with the vector field
 	\begin{align*}
 		\begin{aligned}
 			(x_1\pfrac{b}{x_3}+\bar{u}_{1,1}\pfrac{b}{\bar{u}_1})\partial_{x_2}+(x_1b\pfrac{b}{x_3}+\bar{u}_{1,1}b\pfrac{b}{\bar{u}_1}-1)\partial_{x_3}\,.
 		\end{aligned}
 	\end{align*}
 	Therefore, we obtain the condition
 	\begin{align*}
 		\begin{aligned}
 			\pfrac{b}{\bar{u}_1}(x_1b\pfrac{b}{x_3}\!+\!\bar{u}_{1,1}b\pfrac{b}{\bar{u}_1}\!-\!1)-b\pfrac{b}{\bar{u}_1}(x_1\pfrac{b}{x_3}\!+\!\bar{u}_{1,1}\pfrac{b}{\bar{u}_1})&=0\,,
 		\end{aligned}
 	\end{align*}
 	which simplifies to $\pfrac{b}{\bar{u}_1}=0$. However, $\pfrac{b}{\bar{u}_1}=0$ together with $a=0$ from above is a contradiction to the regularity of the transformation \eqref{eq:inv_trans_1}, which would then read as
	\begin{align*}
	 	\begin{aligned}
		 	u_1&=\bar{u}_2\\
		 	u_2&=b(x)\bar{u}_2\,.
	 	\end{aligned}
 	\end{align*}
 	Hence, the required transformation cannot be of the form \eqref{eq:inv_trans_1}.\\
 	
 	\textit{\textbf{Case 2.}} The second case can be handled analogously. In this case, we have to consider the prolonged system
 	\begin{align}\label{eq:prolonged_2}
 		\begin{aligned}
	 		\dot{x}_1&=h(x,\bar{u}_1,\bar{u}_2)\\
	 		\dot{x}_2&=\bar{u}_2\\
	 		\dot{x}_3&=x_1+\tfrac{\bar{u}_2^2}{2h(x,\bar{u}_1,\bar{u}_2)}\\
	 		\dot{\bar{u}}_1&=\bar{u}_{1,1}\\
	 		\dot{\bar{u}}_{1,1}&=\bar{u}_{1,2}\\
	 		\dot{\bar{u}}_{1,2}&=\bar{u}_{1,3}\,,
 		\end{aligned}
 	\end{align}
	and the corresponding input transformation reads as	
	\begin{align}\label{eq:inv_trans_2}
		\begin{aligned}
			u_1&=h(x,\bar{u}_1,\bar{u}_2)\\
			u_2&=\bar{u}_2\,.
		\end{aligned}
	\end{align}
	By assumption, the prolonged system \eqref{eq:prolonged_2} is static feedback linearizable. Thus, the distributions
	\begin{align*}
		\begin{aligned}
			\mathcal D_p^0&=\Span{\partial_{\bar{u}_{1,3}},\partial_{\bar{u}_2}}\\
			\mathcal D_p^1&=\Span{\partial_{\bar{u}_{1,3}},\partial_{\bar{u}_{1,2}},\partial_{\bar{u}_2},\pfrac{h}{\bar{u}_2}\partial_{x_1}\!\!+\!\partial_{x_2}\!\!+\!\tfrac{2\bar{u}_2h-\bar{u}_2^2\pfrac{h}{\bar{u}_2}}{2h^2}\partial_{x_3}}\\
			&\vdotswithin{=}
		\end{aligned}
	\end{align*}
	are involutive. The involutivity of $\mathcal D_p^1$ implies that
	\begin{align*}
		\begin{aligned}
			\ppfrac{h}{\bar{u}_2}&=0&\text{and}&&\pfrac{}{\bar{u}_2}(\tfrac{2\bar{u}_2h-\bar{u}_2^2\pfrac{h}{\bar{u}_2}}{2h^2})&=0\,.
		\end{aligned}
	\end{align*}
	The first condition implies that $h$ is affine with respect to $\bar{u}_2$, i.e., we actually have $h=a(x,\bar{u}_1)+b(x,\bar{u}_1)\bar{u}_2$. The second condition then yields
	\begin{align*}
		\begin{aligned}
			\pfrac{}{\bar{u}_2}\left(\frac{2a\bar{u}_2+b\bar{u}_2^2}{2(a+b\bar{u}_2)^2}\right)&=0\,,
		\end{aligned}
	\end{align*}
	which simplifies to
	\begin{align*}
		\begin{aligned}
			\frac{a^2}{(a+b\bar{u}_2)^3}&=0\,.
		\end{aligned}
	\end{align*}
	Thus, we must have $a=0$. As a consequence, the distribution $\mathcal D_p^1$ simplifies to
	\begin{align*}
		\begin{aligned}
			\mathcal D_p^1&=\Span{\partial_{\bar{u}_{1,3}},\partial_{\bar{u}_{1,2}},\partial_{\bar{u}_2},b\partial_{x_1}+\partial_{x_2}+\tfrac{1}{2b}\partial_{x_3}}\,,
		\end{aligned}
	\end{align*}
	and for the next distribution we obtain
	\begin{align*}
		\begin{aligned}
			\mathcal D_p^2&=\Span{\partial_{\bar{u}_{1,3}},\partial_{\bar{u}_{1,2}},\partial_{\bar{u}_{1,1}},\partial_{\bar{u}_2},\partial_{\bar{u}_{1,2}},\partial_{\bar{u}_2},\\
				&\hspace{0.4em}b\partial_{x_1}+\partial_{x_2}+\tfrac{1}{2b}\partial_{x_3},\\
				&\hspace{0.4em}(x_1\pfrac{b}{x_3}+\bar{u}_{1,1}\pfrac{b}{\bar{u}_1})\partial_{x_1}-(\tfrac{x_1}{2b^2}\pfrac{b}{x_3}+\tfrac{\bar{u}_{1,1}}{2b^2}\pfrac{b}{\bar{u}_1}+b)\partial_{x_3}}\,,
		\end{aligned}
	\end{align*}
	which must again be involutive. The involutivity implies that
	\begin{align*}
		\pfrac{b}{\bar{u}_1}\partial_{x_1}-\tfrac{1}{2b^2}\pfrac{b}{\bar{u}_1}\partial_{x_3}&\in\mathcal D_p^2\,,
	\end{align*}
	which can only hold if the vector field
	\begin{align*}
		\begin{aligned}
			\pfrac{b}{\bar{u}_1}\partial_{x_1}-\tfrac{1}{2b^2}\pfrac{b}{\bar{u}_1}\partial_{x_3}
		\end{aligned}
	\end{align*}
	is collinear with the vector field
	\begin{align*}
		\begin{aligned}
			(x_1\pfrac{b}{x_3}+\bar{u}_{1,1}\pfrac{b}{\bar{u}_1})\partial_{x_1}-(\tfrac{x_1}{2b^2}\pfrac{b}{x_3}+\tfrac{\bar{u}_{1,1}}{2b^2}\pfrac{b}{\bar{u}_1}+b)\partial_{x_3}\,.
		\end{aligned}
	\end{align*}
	Therefore, we obtain the condition
	\begin{align*}
		\pfrac{b}{\bar{u}_1}(\tfrac{x_1}{2b^2}\pfrac{b}{x_3}\!+\!\tfrac{\bar{u}_{1,1}}{2b^2}\pfrac{b}{\bar{u}_1}\!+\!b)-\tfrac{1}{2b^2}\pfrac{b}{\bar{u}_1}(x_1\pfrac{b}{x_3}\!+\!\bar{u}_{1,1}\pfrac{b}{\bar{u}_1})&=0
	\end{align*}
	which simplifies to $\pfrac{b}{\bar{u}_1}b=0$. Thus, either $\pfrac{b}{\bar{u}_1}=0$ or $b=0$. Since we already have $a=0$, $b=0$ would imply $h=0$, and hence the transformation would not be invertible. However, the other possibility, i.e., $\pfrac{b}{\bar{u}_1}=0$, is in combination with $a=0$ also a contradiction to the regularity of the transformation \eqref{eq:inv_trans_2}, which would then read as
	\begin{align*}
		\begin{aligned}
			u_1&=b(x)\bar{u}_2\\
			u_2&=\bar{u}_2\,.
		\end{aligned}
	\end{align*}
	
	In conclusion, there cannot exist an invertible input transformation \eqref{eq:input_trasnformation_inv} which generates an input such that a three-fold prolongation of this input renders the system \eqref{eq:sys_xuup_intro} static feedback linearizable. The system is thus not linearizable by a three-fold prolongation, and as a consequence of Proposition \ref{prop:linearization_xu} the system \eqref{eq:sys_xuup_intro} cannot be $(x,u)$-flat.\hfill$\square$

\section{CONCLUSIONS}
	We have proven that there exist systems which are flat but do not admit any $(x,u)$-flat output. It thus seems reasonable to conjecture that there also exist systems which are flat but do not admit any $(x,u,\dot{u})$-flat output, and so on. Further research will be devoted to proving this conjecture.
	
\section*{ACKNOWLEDGMENT}
	The authors would like to thank F. Nicolau for interesting discussions on the results presented in this paper.

\addtolength{\textheight}{-2cm}   % This command serves to balance the column lengths

%%%%%%%%%%%%%%%%%%%%%%%%%%%%%%%%%%%%%%%%%%%%%%%%%%%%%%%%%%%%%%%%%%%%%%%%%%%%%%%%

%%%%%%%%%%%%%%%%%%%%%%%%%%%%%%%%%%%%%%%%%%%%%%%%%%%%%%%%%%%%%%%%%%%%%%%%%%%%%%%%

%%%%%%%%%%%%%%%%%%%%%%%%%%%%%%%%%%%%%%%%%%%%%%%%%%%%%%%%%%%%%%%%%%%%%%%%%%%%%%%%
\section*{APPENDIX}

\subsection{Flatness of \eqref{eq:sys_xuup_intro} and Linearizability by Prolongations}
	The flat parameterization of the state $x$ and the input $u$ of the system \eqref{eq:sys_xuup_intro} with respect to the flat output \eqref{eq:xuup_flat_output} is given by
	\begin{align*}
		\begin{aligned}
			x_1&=\tfrac{\dot{y}_1^3\ddot{y}_2 - 3\dot{y}_1^2\ddot{y}_1\dot{y}_2 - \dot{y}_1^2y_1^{(3)}y_2 + 3\dot{y}_1\ddot{y}_1^2y_2 + \ddot{y}_1\ddot{y}_2 - y_1^{(3)}\dot{y}_2}{2(\dot{y}_1^3 + \ddot{y}_1)^2}\\[2ex]
			x_2&=\tfrac{1}{2(\dot{y}_1^3 + \ddot{y}_1)^2}(-\dot{y}_1^4\dot{y}_2 + (y_1\ddot{y}_2 + \ddot{y}_1y_2)\dot{y}_1^3 \\
			&\hspace{2em}- (3\ddot{y}_1\dot{y}_2 + y_1^{(3)}y_2)y_1\dot{y}_1^2 + (3y_1\ddot{y}_1^2y_2 - \ddot{y}_1\dot{y}_2)\dot{y}_1\\
			&\hspace{2em} + y_1\ddot{y}_1\ddot{y}_2 - y_1y_1^{(3)}\dot{y}_2 + \ddot{y}_1^2y_2)\\[2ex]
			x_3&=\tfrac{1}{4(\dot{y}_1^3 + \ddot{y}_1)^2}(2\dot{y}_1^5y_2 - 2y_1\dot{y}_1^4\dot{y}_2 \\
			&\hspace{2em} + (y_1^2\ddot{y}_2 + 2y_1\ddot{y}_1y_2 + 2\dot{y}_2)\dot{y}_1^3\\
			&\hspace{2em} + ((-3y_1^2\dot{y}_2 + 2y_2)\ddot{y}_1 - y_1^2y_2y_1^{(3)})\dot{y}_1^2\\
			&\hspace{2em} + (3y_1^2\ddot{y}_1^2y_2 - 2y_1\ddot{y}_1\dot{y}_2)\dot{y}_1 + 2y_1\ddot{y}_1^2y_2 \\
			&\hspace{2em} + (y_1^2\ddot{y}_2 + 2\dot{y}_2)\ddot{y}_1 - y_1^2y_1^{(3)}\dot{y}_2)\\[1ex]
		\end{aligned}
	\end{align*}
	\begin{align*}
		\begin{aligned}
			u_1 &= \tfrac{1}{2(\dot{y}_1^3 + \ddot{y}_1)^3}(\dot{y}_1^6y_2^{(3)} + (-6\ddot{y}_1\ddot{y}_2 - 4y_1^{(3)}\dot{y}_2 - y_1^{(4)}y_2)\dot{y}_1^5\\
			&\hspace{2em} + (15\ddot{y}_1^2\dot{y}_2 + 10\ddot{y}_1y_1^{(3)}y_2)\dot{y}_1^4 \\
			&\hspace{2em} + (-15\ddot{y}_1^3y_2 + 2\ddot{y}_1y_2^{(3)} - 2y_1^{(3)}\ddot{y}_2 - y_1^{(4)}\dot{y}_2)\dot{y}_1^3 \\
			&\hspace{2em} + (-6\ddot{y}_2\ddot{y}_1^2 + (8y_1^{(3)}\dot{y}_2 - y_1^{(4)}y_2)\ddot{y}_1 + 2y_2(y_1^{(3)})^2)\dot{y}_1^2 \\
			&\hspace{2em} + (-3\ddot{y}_1^3\dot{y}_2 - 2\ddot{y}_1^2y_1^{(3)}y_2)\dot{y}_1 + 3\ddot{y}_1^4y_2 + \ddot{y}_1^2y_2^{(3)} \\ 
			&\hspace{2em} + (-2y_1^{(3)}\ddot{y}_2 - y_1^{(4)}\dot{y}_2)\ddot{y}_1 + 2(y_1^{(3)})^2\dot{y}_2)\\[2ex]
			u_2 &= \tfrac{-y_1}{2(\dot{y}_1^3 + \ddot{y}_1)^3}(-\dot{y}_1^6y_2^{(3)} + (6\ddot{y}_1\ddot{y}_2 + 4y_1^{(3)}\dot{y}_2 + y_1^{(4)}y_2)\dot{y}_1^5 \\
			&\hspace{2em} + (-15\ddot{y}_1^2\dot{y}_2 - 10\ddot{y}_1y_1^{(3)}y_2)\dot{y}_1^4 \\
			&\hspace{2em} + (15\ddot{y}_1^3y_2 - 2\ddot{y}_1y_2^{(3)} + 2y_1^{(3)}\ddot{y}_2 + y_1^{(4)}\dot{y}_2)\dot{y}_1^3 \\
			&\hspace{2em} + (6\ddot{y}_2\ddot{y}_1^2 + (-8y_1^{(3)}\dot{y}_2 + y_1^{(4)}y_2)\ddot{y}_1 - 2y_2(y_1^{(3)})^2)\dot{y}_1^2 \\
			&\hspace{2em} + (3\ddot{y}_1^3\dot{y}_2 + 2\ddot{y}_1^2y_1^{(3)}y_2)\dot{y}_1 - 3\ddot{y}_1^4y_2 - \ddot{y}_1^2y_2^{(3)} \\
			&\hspace{2em} + (2y_1^{(3)}\ddot{y}_2 + y_1^{(4)}\dot{y}_2)\ddot{y}_1 - 2(y_1^{(3)})^2\dot{y}_2)\,.
		\end{aligned}
	\end{align*}
	Above we claimed that the system \eqref{eq:sys_xuup_intro} can be rendered static feedback linearizable by a $4$-fold prolongation of the input $\bar{u}_1=\tfrac{u_2}{u_1}$. Let us explicitly show this. The input $\bar{u}_1$ can be introduced by means of the input transformation
	\begin{align*}
		\begin{aligned}
			\bar{u}_1&=\tfrac{u_2}{u_1}\\
			\bar{u}_2&=u_2
		\end{aligned}
	\end{align*}
	(the particular choice for $\bar{u}_2$ does not matter). Applying this input transformation to \eqref{eq:sys_xuup_intro} and subsequently $4$-fold prolonging $\bar{u}_1$ yields the system
	\begin{align}\label{eq:sys_xuup_intro_prolonged}
		\begin{aligned}
			\dot{x}_1&=\tfrac{\bar{u}_2}{\bar{u}_1}\\
			\dot{x}_2&=\bar{u}_2\\
			\dot{x}_3&=x_1+\tfrac{1}{2}\bar{u}_1\bar{u}_2\,,
		\end{aligned}\qquad
		\begin{aligned}
			\dot{\bar{u}}_1&=\bar{u}_{1,1}\\
			\dot{\bar{u}}_{1,1}&=\bar{u}_{1,2}\\
			\dot{\bar{u}}_{1,2}&=\bar{u}_{1,3}\\
			\dot{\bar{u}}_{1,3}&=\bar{u}_{1,4}
		\end{aligned}
	\end{align}
	with the $n_p=7$-dimensional state $(x_1,x_2,x_3,\bar{u}_1,\bar{u}_{1,1},\bar{u}_{1,2},$ $\bar{u}_{1,3})$ and the input $(\bar{u}_{1,4},\bar{u}_2)$. The distributions involved in the test for static feedback linearizability of \eqref{eq:sys_xuup_intro_prolonged} follow as
	\begin{align*}
		\begin{aligned}
			\mathcal D_p^0&=\Span{\partial_{\bar{u}_{1,4}},\partial_{\bar{u}_2}}\\[2ex]
			\mathcal D_p^1&=\Span{\partial_{\bar{u}_{1,4}},\partial_{\bar{u}_{1,3}},\partial_{\bar{u}_2},\partial_{x_1}+\bar{u}_1\partial_{x_2}+\tfrac{1}{2}\bar{u}_1^2\partial_{x_3}}\\[2ex]
			\mathcal D_p^2&=\Span{\partial_{\bar{u}_{1,4}},\partial_{\bar{u}_{1,3}},\partial_{\bar{u}_{1,2}},\partial_{\bar{u}_2},\partial_{x_1}+\bar{u}_1\partial_{x_2}+\tfrac{1}{2}\bar{u}_1^2\partial_{x_3},\\
			&\hspace{3em}\bar{u}_{1,1}\partial_{x_2}+(\bar{u}_1\bar{u}_{1,1}-1)\partial_{x_3}}\\[2ex]
			\mathcal D_p^3&=\Span{\partial_{\bar{u}_{1,4}},\partial_{\bar{u}_{1,3}},\partial_{\bar{u}_{1,2}},\partial_{\bar{u}_{1,1}},\partial_{\bar{u}_2},\partial_{x_1},\partial_{x_2},\partial_{x_3}}\\[2ex]
			\mathcal D_p^4&=\Span{\partial_{\bar{u}_{1,4}},\partial_{\bar{u}_{1,3}},\partial_{\bar{u}_{1,2}},\partial_{\bar{u}_{1,1}},\partial_{\bar{u}_1},\partial_{\bar{u}_2},\partial_{x_1},\partial_{x_2},\partial_{x_3}}\,.
		\end{aligned}
	\end{align*}
	It is easy to see that all these distributions are involutive and we have $\Dim{\mathcal D_p^4}=n_p+2$. Thus, according to Theorem \ref{thm:sfl}, the prolonged system \eqref{eq:sys_xuup_intro_prolonged} is static feedback linearizable. A possible linearizing output of \eqref{eq:sys_xuup_intro_prolonged} follows as 
	\begin{align*}
		\begin{aligned}
			y_1&=\bar{u}_1\\
			y_2&=(\bar{u}_1^2x_1 - 2\bar{u}_1x_2 + 2x_3)\bar{u}_{1,1} - 2\bar{u}_1x_1 + 2x_2\,,
		\end{aligned}
	\end{align*}
	which is exactly the flat output \eqref{eq:xuup_flat_output}. Indeed, substituting $\bar{u}_1=\tfrac{u_2}{u_1}$ and $\bar{u}_{1,1}=\dot{\bar{u}}_1=\tfrac{1}{u_1^2}(\dot{u}_2u_1-u_2\dot{u}_1)$ into the linearizing output yields \eqref{eq:xuup_flat_output}.

\bibliographystyle{IEEEtran} 
\bibliography{IEEEabrv,Bibliography}

\end{document}